\documentclass[9pt]{article}
\usepackage{amsmath}
\usepackage[dvips]{graphicx}
\usepackage{textcomp}
\usepackage{amsbsy}
\usepackage{latexsym}
\usepackage[mathscr]{eucal}
\usepackage{amsfonts}
\usepackage{amssymb}
\usepackage[usenames]{color}
\usepackage{stmaryrd}

\usepackage{fullpage}

\begin{document}
\bibliographystyle{plain}
\title
{
On the stability and accuracy of least squares approximations
}
\author{
Albert Cohen, Mark A. Davenport and Dany Leviatan\thanks{%
This research has been partially supported by the
ANR Defi08 ``ECHANGE'' and by the US NSF grant DMS-1004718. Portions
of this work were completed while M.A.D.\ and D.L.\ were visitors at
University Pierre et Marie Curie} }
\hbadness=10000
\vbadness=10000
\newtheorem{lemma}{Lemma}
\newtheorem{prop}[lemma]{Proposition}
\newtheorem{cor}[lemma]{Corollary}
\newtheorem{theorem}{Theorem}
\newtheorem{remark}[lemma]{Remark}
\newtheorem{example}[lemma]{Example}
\newtheorem{definition}[lemma]{Definition}
\newtheorem{proper}[lemma]{Properties}
\newtheorem{assumption}[lemma]{Assumption}
%
\def\RR{\rm \hbox{I\kern-.2em\hbox{R}}}
\def\NN{\rm \hbox{I\kern-.2em\hbox{N}}}
\def\ZZ{\rm {{\rm Z}\kern-.28em{\rm Z}}}
\def\CC{\rm \hbox{C\kern -.5em {\raise .32ex \hbox{$\scriptscriptstyle
|$}}\kern
-.22em{\raise .6ex \hbox{$\scriptscriptstyle |$}}\kern .4em}}
\def\vp{\varphi}
\def\<{\langle}
\def\>{\rangle}
\def\t{\tilde}
\def\i{\infty}
\def\e{\varepsilon}
\def\sm{\setminus}
\def\nl{\newline}
\def\o{\bar}
\def\wt{\widetilde}
\def\wh{\widehat}
\def\cT{{\cal T}}
\def\cA{{\cal A}}
\def\cI{{\cal I}}
\def\cV{{\cal V}}
\def\cB{{\cal B}}
\def\cR{{\cal R}}
\def\cD{{\cal D}}
\def\cP{{\cal P}}
\def\cJ{{\cal J}}
\def\cM{{\cal M}}
\def\cO{{\cal O}}
\def\Chi{\raise .3ex
\hbox{\large $\chi$}} \def\vp{\varphi}
\def\lsima{\hbox{\kern -.6em\raisebox{-1ex}{$~\stackrel{\textstyle<}{\sim}~$}}\kern -.4em}
\def\lsim{\hbox{\kern -.2em\raisebox{-1ex}{$~\stackrel{\textstyle<}{\sim}~$}}\kern -.2em}
\def\[{\Bigl [}
\def\]{\Bigr ]}
\def\({\Bigl (}
\def\){\Bigr )}
\def\[{\Bigl [}
\def\]{\Bigr ]}
\def\({\Bigl (}
\def\){\Bigr )}
\def\L{\pounds}
\def\pr{{\rm Prob}}

\newcommand{\cs}[1]{{\color{blue}{#1}}}
\def\ds{\displaystyle}
\def\ev#1{\vec{#1}}     
\newcommand{\lt}{\ell^{2}(\nabla)}
\def\Supp#1{{\rm supp\,}{#1}}
\def\R{\mathbb{R}}
\def\P{\mathbb{P}}
\def\E{\mathbb{E}}
\def\nl{\newline}
\def\T{{\relax\ifmmode I\!\!\hspace{-1pt}T\else$I\!\!\hspace{-1pt}T$\fi}}
\def\N{\mathbb{N}}
\def\Z{\mathbb{Z}}
\def\N{\mathbb{N}}
\def\Zd{\Z^d}
\def\Q{\mathbb{Q}}
\def\C{\mathbb{C}}
\def\Rd{\R^d}
\def\gsim{\mathrel{\raisebox{-4pt}{$\stackrel{\textstyle>}{\sim}$}}}
\def\sime{\raisebox{0ex}{$~\stackrel{\textstyle\sim}{=}~$}}
\def\lsim{\raisebox{-1ex}{$~\stackrel{\textstyle<}{\sim}~$}}
\def\div{\mbox{ div }}
\def\M{M}  \def\NN{N}                  
\def\L{{\ell}}               
\def\Le{{\ell^1}}            
\def\Lz{{\ell^2}}
\def\Let{{\tilde\ell^1}}     
\def\Lzt{{\tilde\ell^2}}
\def\Ltw{\ell^\tau^w(\nabla)}
\def\t#1{\tilde{#1}}
\def\la{\lambda}
\def\La{\Lambda}
\def\ga{\gamma}
\def\BV{{\rm BV}}
\def\Ga{\eta}
\def\al{\alpha}
\def\cZ{{\cal Z}}
\def\cA{{\cal A}}
\def\cU{{\cal U}}
\def\cV{{\cal V}}
\def\argmin{\mathop{\mathrm argmin}}
\def\argmax{\mathop{\mathrm argmax}}
\def\prob{\mathop{\rm prob}}

\def\cO{{\cal O}}
\def\cA{{\cal A}}
\def\cC{{\cal C}}
\def\cF{{\cal F}}
\def\bu{{\bf u}}
\def\bz{{\bf z}}
\def\bZ{{\bf Z}}
\def\bI{{\bf I}}
\def\cE{{\cal E}}
\def\cD{{\cal D}}
\def\cG{{\cal G}}
\def\cI{{\cal I}}
\def\cJ{{\cal J}}
\def\cM{{\cal M}}
\def\cN{{\cal N}}
\def\cT{{\cal T}}
\def\cU{{\cal U}}
\def\cV{{\cal V}}
\def\cW{{\cal W}}
\def\cL{{\cal L}}
\def\cB{{\cal B}}
\def\cG{{\cal G}}
\def\cK{{\cal K}}
\def\cS{{\cal S}}
\def\cP{{\cal P}}
\def\cQ{{\cal Q}}
\def\cR{{\cal R}}
\def\cU{{\cal U}}
\def\bL{{\bf L}}
\def\bl{{\bf l}}
\def\bK{{\bf K}}
\def\bC{{\bf C}}
\def\X{X\in\{L,R\}}
\def\ph{{\varphi}}
\def\D{{\Delta}}
\def\H{{\cal H}}
\def\bM{{\bf M}}
\def\bx{{\bf x}}
\def\bj{{\bf j}}
\def\bG{{\bf G}}
\def\bX{{\bf X}}
\def\bP{{\bf P}}
\def\bW{{\bf W}}
\def\bT{{\bf T}}
\def\bV{{\bf V}}
\def\bv{{\bf v}}
\def\bn{{\bf n}}
\def\bt{{\bf t}}
\def\bz{{\bf z}}
\def\bw{{\bf w}}
\def \meas {{\rm meas}}
\def\rhom{{\rho^m}}
\def\lll{\langle}
\def\argmin{\mathop{\rm argmin}}
\def\argmax{\mathop{\rm argmax}}
\def\dJ{\nabla}
\newcommand{\ba}{{\bf a}}
\newcommand{\bb}{{\bf b}}
\newcommand{\bc}{{\bf c}}
\newcommand{\bd}{{\bf d}}
\newcommand{\bs}{{\bf s}}
\newcommand{\bff}{{\bf f}}
\newcommand{\bp}{{\bf p}}
\newcommand{\bg}{{\bf g}}
\newcommand{\by}{{\bf y}}
\newcommand{\br}{{\bf r}}
\newcommand{\be}{\begin{equation}}
\newcommand{\ee}{\end{equation}}
\newcommand{\bea}{$$ \begin{array}{lll}}
\newcommand{\eea}{\end{array} $$}
\def \Vol{\mathop{\rm  Vol}}
\def \mes{\mathop{\rm mes}}
\def \Prob{\mathop{\rm  Prob}}
\def \exp{\mathop{\rm    exp}}
\def \sign{\mathop{\rm   sign}}
\def \sp{\mathop{\rm   span}}
\def \vphi{{\varphi}}
\def \csp{\overline \mathop{\rm   span}}
\newcommand{\KL}{Karh\'unen-Lo\`eve }
%
\newcommand{\beqn}{\begin{equation}}
\newcommand{\eeqn}{\end{equation}}
\def\beginproof{\noindent{\bf Proof:}~ }
\def\endproof{\hfill\rule{1.5mm}{1.5mm}\\[2mm]}

\newenvironment{Proof}{\noindent{\bf Proof:}\quad}{\endproof}

\renewcommand{\theequation}{\thesection.\arabic{equation}}
\renewcommand{\thefigure}{\thesection.\arabic{figure}}

\makeatletter
\@addtoreset{equation}{section}
\makeatother

\newcommand\abs[1]{\left|#1\right|}
\newcommand\clos{\mathop{\rm clos}\nolimits}
\newcommand\trunc{\mathop{\rm trunc}\nolimits}
\renewcommand\d{d}
\newcommand\dd{d}
\newcommand\diag{\mathop{\rm diag}}
\newcommand\dist{\mathop{\rm dist}}
\newcommand\diam{\mathop{\rm diam}}
\newcommand\cond{\mathop{\rm cond}\nolimits}
\newcommand\eref[1]{{\rm (\ref{#1})}}
\newcommand{\iref}[1]{{\rm (\ref{#1})}}
\newcommand\Hnorm[1]{\norm{#1}_{H^s([0,1])}}
\def\int{\intop\limits}
\renewcommand\labelenumi{(\roman{enumi})}
\newcommand\lnorm[1]{\norm{#1}_{\ell^2(\Z)}}
\newcommand\Lnorm[1]{\norm{#1}_{L^2([0,1])}}
\newcommand\LR{{L^2(\R)}}
\newcommand\LRnorm[1]{\norm{#1}_\LR}
\newcommand\Matrix[2]{\hphantom{#1}_#2#1}
\newcommand\norm[1]{\left\|#1\right\|}
\newcommand\ogauss[1]{\left\lceil#1\right\rceil}
\newcommand{\QED}{\hfill
\raisebox{-2pt}{\rule{5.6pt}{8pt}\rule{4pt}{0pt}}%
  \smallskip\par}
\newcommand\Rscalar[1]{\scalar{#1}_\R}
\newcommand\scalar[1]{\left(#1\right)}
\newcommand\Scalar[1]{\scalar{#1}_{[0,1]}}
\newcommand\supp{\mathop{\rm supp}}
\newcommand\ugauss[1]{\left\lfloor#1\right\rfloor}
\newcommand\with{\, : \,}
\newcommand\Null{{\bf 0}}
\newcommand\bA{{\bf A}}
\newcommand\bB{{\bf B}}
\newcommand\bR{{\bf R}}
\newcommand\bD{{\bf D}}
\newcommand\bE{{\bf E}}
\newcommand\bF{{\bf F}}
\newcommand\bH{{\bf H}}
\newcommand\bU{{\bf U}}
\newcommand\cH{{\cal H}}
\newcommand\sinc{{\rm sinc}}
\def\enorm#1{| \! | \! | #1 | \! | \! |}

\newcommand{\dm}{\frac{d-1}{d}}

\let\bm\bf
\newcommand{\bbeta}{{\mbox{\boldmath$\beta$}}}
\newcommand{\bal}{{\mbox{\boldmath$\alpha$}}}
\newcommand{\bbi}{{\bm i}}

\def\nnew{\color{Red}}
\def\mnew{\color{Blue}}

\newcommand{\dI}{\Delta}

\maketitle
\date{}

\begin{abstract}
We consider the problem of reconstructing an unknown function $f$ on
a domain $X$ from samples of $f$ at $n$ randomly chosen
points with respect to a given measure $\rho_X$. Given a sequence of
linear spaces $(V_m)_{m>0}$ with ${\rm dim}(V_m)=m\leq n$, we study the
least squares approximations from the spaces $V_m$. It is well known
that such approximations can be inaccurate when $m$ is too close to
$n$, even when the samples are noiseless. Our main result provides
a criterion on $m$ that describes the needed amount of
regularization to ensure that the least squares method is stable
and that its accuracy, measured in $L^2(X,\rho_X)$, is comparable to
the best approximation error of $f$ by elements from $V_m$.  We illustrate this
criterion for various approximation schemes, such as trigonometric
polynomials, with $\rho_X$ being the uniform measure, and algebraic
polynomials, with $\rho_X$ being either the uniform or Chebyshev
measure. For such examples we also prove similar stability results
using deterministic samples that are equispaced with respect to these measures.
\end{abstract}


\section{Introduction and main results}

Let $X$ be a domain of $\R^d$ and $\rho_X$ be a probability measure
on $X$. We consider the problem of estimating an unknown function $f: X \to \R$
from samples $(y_i)_{i=1,\dots,n}$ which are either noiseless or noisy observations
of $f$ at the points $(x_i)_{i=1,\dots,n}$, where the $x_i$ are i.i.d.\ with respect to $\rho_X$.  We measure
the error between $f$ and its estimator $\t f$ in the $L^2(X,\rho_X)$ norm
$$
\|v \|:=\(\int_X |v(x)|^2 d\rho_X(x)\)^{1/2},
$$
and we denote by $\<\cdot,\cdot\>$ the associated inner product.

Given a fixed sequence of finite dimensional spaces $(V_m)_{m\geq 1}$ of $L^2(X,\rho_X)$ such that $\dim(V_m)=m$.
We would like to compute the best approximation of $f$ in $V_m$.  This is given by the $L^2(X,\rho_X)$ orthogonal projector onto $V_m$, which we denote by $P_m$:
$$
P_mf:= \argmin_{v\in V_m}\|f-v\|.
$$
We let
$$
e_m(f)=\|f-P_mf\|
$$
denote the best approximation error.

In general, we may not have access to either $\rho_X$ or any information about $f$ aside from the observations at the points $(x_i)_{i=1,\dots,n}$.  In this case we cannot explicitly compute $P_m f$.  A natural approach in this setting is to consider the solution of the least squares problem
$$
w = \argmin_{v\in V_m} \sum_{i=1}^n |y_i-v(x_i)|^2.
$$
Typically, we are interested in the case where $m \leq n$ which is the regime where this problem may admit a unique solution.

In the noiseless case $y_i=f(x_i)$, and hence $w$ may be viewed as the
application of the least squares projection operator onto $V_m$ to $f$, i.e., we can write
$$
w = P^n_m f:= \argmin_{v\in V_m}\|f-v\|_n
$$
where
$$
\|v \|_n:=\(\frac 1 n\sum_{i=1}^n |v(x_i)|^2\)^{1/2}
$$
is the $L^2$ norm with respect to the empirical measure and,
analogously, $\<\cdot,\cdot\>_n$ the associated empirical inner
product.

It is well known that least squares approximations may be inaccurate
even when the measured samples are noiseless. For example, if $V_m$
is the space $\P_{m-1}$ of algebraic polynomials of degree $m-1$
over the interval $[-1,1]$ and if we choose $m=n$, this corresponds
to Lagrange interpolation, which is known to be highly unstable, failing
to converge towards $f$ when given values at uniformly spaced samples, even when
$f$ is infinitely smooth (the ``Runge phenomenon''). Regularization by taking $m$ substantially smaller than $n$ may therefore be needed even in a noise-free
context. The goal of this paper is to provide a mathematical analysis
on the exact needed amount of such regularization.
\nl
\nl
{\bf Stability of the least squares problem.} The solution of the least squares problem can be computed by solving an
$m\times m$ system: specifically, if $(L_1,\dots,L_m)$ is an arbitrary basis for
$V_m$, then we can write
$$
w=\sum_{j=1}^m u_j L_j,
$$
where $\bu=(u_j)_{j=1,\dots,m}$ is the solution of the $m\times m$
system \be \label{system} \bG\bu =\bff, \;\; \ee with
$\bG:=(\<L_j,L_k\>_n)_{j,k=1,\dots,m}$ and $\bff= (\frac 1
n\sum_{i=1}^n y_i L_k(x_i))_{k=1,\dots, m}$. In the noiseless case
$y_i=f(x_i)$, so that we can also write
$\bff:=(\<f,L_k\>_n)_{k=1,\dots,m}$. In the event that $\bG$ is
singular, we simply set $w=0$.

For the purposes of our analysis, suppose that the basis
$(L_1,\dots,L_m)$ is orthonormal in the sense of $L^2(X,\rho_X)$.\footnote{While such a basis is generally not accessible when $\rho_X$ is
unknown, we require it only for the analysis. The actual computation
of the estimator can be made using any known basis of $V_m$, since
the solution $w$ is independent of the basis used in computing it.} In
this case we have
$$
\E(\bG)=(\<L_j,L_k\>)_{j,k=1,\dots,m} = \bI.
$$
Our analysis requires an understanding of how the random matrix $\bG$ deviates from
its expectation $\bI$ in probability. Towards this end, we introduce the quantity
$$
K(m):=\sup_{x\in X}\sum_{j=1}^m |L_j(x)|^2.
$$
Note that the function $\sum_{j=1}^m |L_j(x)|^2$ is invariant with respect to a rotation
applied to $(L_1,\dots,L_m)$ and therefore independent of the choice
of the orthonormal basis: it only depends on the space $V_m$
and on the measure $\rho_X$, and hence $K(m)$ also depends only on $V_m$ and $\rho_X$.
Also note that
$$
K(m)\geq \sum_{j=1}^m \|L_j\|^2= m.
$$
We also will use the notation
$$
\interleave \bM \interleave=\max_{\bv\neq 0} \frac {|\bM \bv|} {|\bv|},
$$
for the spectral norm of a matrix.

Our first result is a probabilistic estimate of the comparability of the norms
$\|\cdot\|$ and $\|\cdot\|_n$ uniformly over the space $V_m$.
This is equivalent to the proximity of the matrices $\bG$ and $\bI$ in spectral norm,
since we have that for all $\delta\in [0,1]$,
$$
\interleave\bG-\bI\interleave\leq \deltaÊ\Leftrightarrow \left
|\|v\|_n^2-\|v\|^2\right | \leq \deltaÊ\|v\|^2, \;\; v\in V_m.
$$
\begin{theorem}
\label{theomat} For $0<\delta<1$, one has the estimate \be {\rm
Pr}\,\{\interleave\bG-\bI\interleave > \delta\}={\rm Pr}\,\{\exists v\in V_m\; : \;
\left |\|v\|_n^2-\|v\|^2\right | > \delta \|v\|^2\} \leq 2m
\exp\left \{-\frac {c_\delta n}{K(m)}\right \}, \label{taildelta}
\ee where $c_\delta:=(1+\delta)\log(1+\delta)-\delta>0$.
\end{theorem}

The proof of Theorem \ref{theomat} is a simple application of tail
bounds for sums of random matrices obtained in \cite{AW}.  A
consequence of this result is that the norms $\|\cdot\|$ and
$\|\cdot\|_n$ are comparable with high probability if $K(m)$ is
smaller than $n$ by a logarithmic factor: for example taking
$\delta=\frac 12$, we find that for any $r>0$,
\be
{\rm Pr} \, \left\{ \interleave\bG-\bI\interleave > \frac 1 2 \right\}={\rm Pr}\, \left\{\exists v\in V_m\; :
\; \left |\|v\|_n^2-\|v\|^2\right | > \frac 1 2 \|v\|^2 \right\} \leq
2n^{-r}, \label{tailhalf}
\ee
if $m$ is such that
\be
K(m) \leq \kappa\frac n {\log n},\;\; {\rm with}\;\; \kappa:=\frac
{c_{1/2}}{1+r}= \frac {3\log(3/2)-1}{2+2r}. \label{condm}
\ee

The above condition thus ensures that $\bG$ is well conditioned with high probability.
It can also be thought of as ensuring that
the least squares problem is {\it stable} with high probability.
Indeed the right side of the least squares system can be written as $\bff=\bM\by$ with
$$
\bM=\frac 1 n(L_j(x_i))_{j,i\in \{1,\dots,m\}\times \{1,\dots,n\}},
$$
an $m\times n$ matrix. Observing that
$\<\bG\bv,\bv\>=n|\bM^T\bv|^2$, we find that
$$
\interleave\bM\interleave=\interleave\bM^T\interleave=\(\frac 1  {n} \interleave\bG\interleave\)^{1/2}.
$$
Therefore, if $\interleave\bG-\bI\interleave \leq \frac 1 2$, then we have that for any
data vector $\by$ the solution $w=\sum_{j=1}^m u_j L_j$ satisfies
$$
\|w\|=|\bu|\leq \interleave\bG^{-1}\interleave\cdot\interleave\bM\interleave\cdot|\by| \leq \frac 1
{\sqrt n}2\sqrt{\frac 3 2}|\by|,
$$
which thus gives the stability estimate
$$
\|w\|\leq C\(\frac 1 n \sum_{i=1}^n |y_i|^2\)^{1/2},\;\; C=\sqrt 6.
$$
In the noiseless case, this can be written as $\|P^n_mf\| \leq C\|f\|_n$, i.e., the least squares projection
is stable between the norms $\|\cdot\|_n$ and $\|\cdot\|$. Note that since $K(m)$ not only depends
on $V_m$ but also on the measure $\rho_X$, the range of $m$ such that
the condition \iref{condm} holds is strongly tied to the choice of the measure.
This issue is illustrated further in our numerical experiments.

Let us mention that similar probabilistic bounds have been previously obtained,
see in particular \S 5.2 in \cite{B}. These earlier results allow
us to obtain the bound \iref{tailhalf}, however relying on
the stronger condition
$$
K(m)\lsim \(\frac n {\log (n)}\)^{1/2}.
$$
The numerical results for polynomial least squares that we present in \S 3 hint that the weaker condition $K(m)\lsim \frac n {\log (n)}$ is sharp.  The quantity $K(m)$ was also used in \cite{BM} in order to control the $L^\infty(X)$ norm
and the $L^2(\rho_X)$ norm.
\nl
\nl
{\bf Accuracy of least squares approximation.}  As an application, we can derive an estimate for the
error of least squares approximation in expectation.
Here, we make the assumption that a
uniform bound
\be
|f(x)| \leq L,
\label{unifbound}
\ee
holds for almost every $x$ with respect to $\rho_X$. For $m\leq n$, we consider
the truncated least squares estimator
$$
\t f=T_L(w),
$$
where $T_L(t)={\rm sign}(t)\max\{L,|t|\}$. Our first result deals with the noiseless case.

\begin{theorem}
\label{theonoiseless}
In the noiseless case, for any $r>0$, if $m$ is such that the condition \iref{condm} holds, then
\be
\E(\|f-\t f\|^2)\leq (1+\e(n))e_m(f)^2+8L^2n^{-r},
\ee
where $\e(n):= \frac {4\kappa} {\log (n)}\to 0$ as $n\to +\infty$, with $\kappa$ as in \iref{condm}
\end{theorem}

At this point a few remarks are due regarding
the implications of this result in terms of the convergence
rate of the estimate.

Consider the following general setting of regression
on a random design: we observe independent
samples
\be
z_i=(x_i,y_i)_{i=1,\dots,n}
\ee
of a variable $z=(x,y)$ of law $\rho$ over $X\times Y$
and marginal law $\rho_X$ over $X$,
and we want to estimate from these samples
the regression function
defined as the conditional expectation
\be
f(x):=\E(y|x).
\ee
We assume that the maximal variance
\be
\label{maxvar}
\sigma^2:=\sup_{x\in X} \E (\left. |y-f(x)|^2  \right | x ),
\ee
is bounded. We thus think of the $y_i$ as noisy observations
of $f$ at $x_i$ with additive noise of variance at most $\sigma^2$,
namely
\be
y_i=f(x_i)+\eta_i,
\ee
where the $\eta_i$ are independent realizations of the variable $\eta:=y-f(x)$.

Assuming that $f$ satisfies the uniform bound \iref{unifbound},
one computes the truncated least squares
estimator now with $y_i$ in place of $f(x_i)$. A typical
convergence bound for this estimator, see for example Theorem 11.3 in \cite{GKKW}, is
\be
\label{typicalbound}
\E(\|f-\t f\|^2)\leq C\(e_m(f)^2+\max\{L^2,\sigma^2\}\frac {m\log n}{n}\).
\ee
Convergence rates may be found after balancing the two terms, but
they are limited by the optimal learning rate $n^{-1}$, and
this limitation persists even in the noiseless case $\sigma^2=0$
due to the presence of $L^2$ in the right side of \iref{typicalbound}.
In contrast, Theorem \ref{theonoiseless} yields fast convergence rates,
provided that the approximation error $e_m$ has fast decay and that the value of $m$
satisfying \iref{condm} can be chosen large enough.

One motivation for studying
the noiseless case is the numerical treatment of
parameter dependent PDEs of the general form
$$
\cF(f,x)=0,
$$
where $x$ is a vector of parameters in some compact set
$\cP\in\R^d$. We can consider the solution map $x\mapsto f(x)$ either as
giving the exact solution to the PDE for the given value of the
parameter vector $x$ or as the exact result of a numerical solver
for this value of $x$. In the stochastic PDE context, $x$ is random
and obeys a certain law which may be known or unknown. From a random
draw $(x_i)_{i=1,\dots,n}$, we obtain solutions $f_i=f(x_i)$ which
are noiseless observations of the solution map, and are interested
in reconstructing this map. In instances such as elliptic problems
with parameters in the diffusion coefficients, the solution map can
be well-approximated by polynomials in $x$ (see \cite{CDS}). In this
context, an initial study of the needed amount of regularization was
given in \cite{MNST}, however specifically targeted towards
polynomial least squares.

For the noisy regression problem described above, our
analysis can also be adapted in order to derive the following result.

\begin{theorem}
\label{theonoisy} For any $r>0$, if $m$ is such that the condition
\iref{condm} holds, then
\be \E(\|f-\t f\|^2)\leq
(1+2\e(n))e_m(f)^2+8L^2n^{-r}+8\sigma^2 \frac m n,
\ee
with $\e(n)$
as in Theorem $\ref{theonoiseless}$ and $\sigma$ is the maximal
variance given by \iref{maxvar}.
\end{theorem}

In the noiseless case, the bound in Theorem \ref{theonoiseless}
suggests that $m$ should be chosen as large as possible under the
constraint that \iref{condm} holds. In the noisy case, the value of
$m$ minimizing the bound in Theorem \ref{theonoisy} also depends on
the decay of $e_m$, which is generally unknown. In such a situation,
a classical way of choosing the value of $m$ is by a model selection
procedure, such as adding a complexity penalty in the least squares
or using an independent validation sample. Such procedures
can also be of interest in the noiseless case when the measure
$\rho_X$ is unknown, since the maximal value of $m$ such
that \iref{condm} holds is then also unknown.

Let us give an example of how the results in Theorems \ref{theonoiseless}
and \ref{theonoisy} lead to specific rates of convergence
in terms of the number of samples: assume that $X=[-1,1]$ is
equipped with the uniform measure $\rho_X=\frac {dx} 2$ and
that $V_m=\P_{m-1}$ is the space of algebraic polynomials of
degree $m-1$. Then, if $f$ belongs to $C^r(X)$ the space of $r$-times differentiable functions,
it is well-known that $e_m(f)^2 \lsim m^{-2r}$. On the one hand
the results in \S 3 show that condition \iref{condm}
can be ensured with $m\sim (n/\log n)^{1/2}$.
Therefore, in the noiseless case, we obtain a bound proportional to
$n^{-r}$ for the mean squared error, up the logarithmic factor.
In the noisy case, after balancing the approximation and variance terms,
we obtain a bound proportional to $\sigma^{2r/(r+1)}n^{-r/(r+1)}$.
On the other hand, these rates can be improved
with $r$ replaced by $2r$ if we use
the Chebyshev non-uniform measure that concentrates near the end-points,
since in that case the results in \S 3 show that
condition \iref{condm} can be ensured with $m\sim n/\log n$.

The rest of our paper is organized as follows: we give the
proofs of the above results in \S 2 and we present
in \S 3 examples of applications to classical
approximation schemes such as piecewise constants,
trigonometric polynomials, or algebraic polynomials.
For such examples, we study the range of $m$ such that
\iref{condm} holds and show that this range is in accordance
with stability results that can be proved for deterministic sampling.
Numerical illustrations are given for algebraic polynomial approximation.

\section{Proofs}

{\bf Proof of Theorem \ref{theomat}:} The matrix $\bG$ can be written
as
$$
\bG=\bX_1+\cdots+\bX_n,
$$
where the $\bX_i$ are i.i.d. copies of the random matrix
$$
\bX=\frac 1 n (L_j(x)L_k(x))_{j,k=1,\dots,m},
$$
where $x$ is distributed according to $\rho_X$.
We use the following Chernoff bound from \cite{T},
originally obtained by \cite{AW}: if
$\bX_1,\dots,\bX_n$ are independent $m\times m$
random self-adjoint and positive matrices satisfying
$$
\lambda_{\max}(\bX_i)=\interleave\bX_i\interleave \leq R,
$$
almost surely,
then with
$$
\mu_{\min}:=\lambda_{\min}\(\sum_{i=1}^n \E(\bX_i)\)\quad \quad {\rm and}
\quad \quad \mu_{\max}:=\lambda_{\max}\(\sum_{i=1}^n \E(\bX_i)\),
$$
one has
$$
{\rm Pr}\,\left \{ \lambda_{\min}\(\sum_{i=1}^n \bX_i\)\leq (1-\delta)\mu_{\min}\right \} \leq m\(\frac {e^{-\delta}}{(1-\delta)^{1-\delta}}\)^{\mu_{\min}/R},\;\; 0\leq \delta<1,
$$
and
$$
{\rm Pr}\,\left \{ \lambda_{\max}\(\sum_{i=1}^n \bX_i\)\geq (1+\delta)\mu_{\max}\right\} \leq m\(\frac {e^{\delta}}{(1+\delta)^{1+\delta}}\)^{\mu_{\max}/R},\;\; \delta \geq 0
$$
In our present case, we have $\sum_{i=1}^n \E(\bX_i)=n\E(\bX)=\bI$
so that $\mu_{\min}=\mu_{\max}=1$. It is easily checked that $\frac
{e^{\delta}}{(1+\delta)^{1+\delta}}\geq \frac
{e^{-\delta}}{(1-\delta)^{1-\delta}}$ for $0<\delta<1$, and
therefore
$$
{\rm Pr}\, \{\interleave\bG-\bI\interleave > \delta\} \leq 2 m \(\frac
{e^{\delta}}{(1+\delta)^{1+\delta}}\)^{1/R} =2m \exp\(-\frac
{c_\delta} R\).
$$
We next use the fact that a rank 1 symmetric matrix $ab^T=(b_ja_k)_{j,k=1,\dots,m}$ has
its spectral norm equal to the product of the Euclidean norms of the vectors $a$ and $b$,
and therefore
$$
\interleave\bX\interleave \leq \frac 1 n \sum_{j=1}^m |L_j(x)|^2=\frac{K(m)} n,
$$
almost surely. We may therefore take $R=\frac{K(m)} n$ which concludes the proof. \hfill $\Box$
\nl
\nl
{\bf Proof of Theorem \ref{theonoiseless}:}
We denote by $d\rho_X^n:=\otimes^n d\rho_X$ the probability
measure of the draw. We also denote by $\Omega$ the set of all possible draws, that we divide into
the set $\Omega_+$ of all draw such that
$$
\interleave\bG-\bI\interleave \leq \frac 1 2,
$$
and the complement set
$\Omega_-:=\Omega\sm\Omega_+$.
According to (\ref{tailhalf}),
we have
\be
{\rm Pr}\{\Omega_-\}=\int_{\Omega_-} d\rho_X^n\leq 2n^{-r},
\label{probomegaminus}
\ee
under the condition \iref{condm}. This leads to
$$
\E(\|f-\t f\|^2)=\int_\Omega \|f-\t f\|^2d\rho_X^n
\leq \int_{\Omega_+} \|f-P_m^nf\|^2d\rho_X^n+ 8L^2n^{-r},
$$
where we have used $\|f-\t f\|^2 \leq 2L^2$,
as well as the fact that $T_L$ is a contraction that preserves $f$.

It remains to prove that the first term in the above right side
is bounded by $(1+\e(n))e_m(f)^2$. With $g:=f-P_mf$,
we observe that
$$
f-P_m^nf=f-P_m f+P_m^nP_m f-P_m^nf=g-P_m^n g.
$$
Since $g$ is orthogonal to $V_m$, we thus have
$$
\|f-P_m^nf\|^2 =\|g\|^2+\|P_m^n g\|^2=\|g\|^2+\sum_{j=1}^m |a_j|^2,
$$
where
$\ba=(a_j)_{j=1,\dots,m}$ is solution of the system
$$
\bG\ba=\bb, \;\;
$$
with $\bb:=(\<g,L_k\>_n)_{k=1,\dots,m}$. When the draw belongs to $\Omega_+$, we have
$\|\bG^{-1}\|_2\leq 2$ and therefore
$$
\sum_{j=1}^m |a_j|^2 \leq 4 \sum_{k=1}^m |\<g,L_k\>_n|^2.
$$
It follows that
$$
\int_{\Omega_+} \|f-P_m^nf\|^2d\rho_X^n\leq
\int_{\Omega_+} \left( \|g\|^2+4 \sum_{k=1}^m |\<g,L_k\>_n|^2 \right)d\rho_X^n
\leq \|g\|^2+4 \sum_{k=1}^m\E( |\<g,L_k\>_n|^2).
$$
We estimate each of the $\E( |\<g,L_k\>_n|^2)$ as follows:
\begin{align*}
\E( |\<g,L_k\>_n|^2) & =\frac 1 {n^2}\sum_{i=1}^n\sum_{j=1}^n \E(g(x_i)g(x_j)L_k(x_i)L_k(x_j)) \\
 & =\frac 1 {n^2} \(n(n-1)|\E(g(x)L_k(x))|^2+n\E(|g(x)L_k(x)|^2)\) \\
 & = \(1-\frac 1 n\)|\<g,L_k\>|^2 +\frac 1 n\int_X |g(x)|^2|L_k(x)|^2d\rho_X \\
 &=\frac 1 n\int_X |g(x)|^2|L_k(x)|^2d\rho_X,
\end{align*}
where we have used the fact that $g$ is orthogonal to $V_m$ and thus to $L_k$.
Summing over $k$, we obtain
$$
\sum_{k=1}^m\E( |\<g,L_k\>_n|^2)\leq \frac {K(m)} n \|g\|^2
\leq \frac \kappa {\log (n)}\|g\|^2,
$$
where we have used \iref{condm}. We have thus proven that
$$
\int_{\Omega_+} \|f-P_m^nf\|^2d\rho_X^n\leq (1+\frac {4\kappa} {\log (n)})\|g\|^2=(1+\e (n))e_m(f)^2,
$$
which concludes the proof. \hfill $\Box$
\nl
\nl
{\bf Proof of Theorem \ref{theonoisy}:} We define the additive noise in the sample by writing
$$
y_i=f(x_i)+\eta_i,
$$
and thus the $\eta_i$ are i.i.d.\ copies of the variable
$$
\eta=y-f(x).
$$
Note that $\eta$ and $x$ are not assumed to be independent. However
we have
$$
\E(\eta |x)=0,
$$
which implies the decorrelation property
$$
\E(\eta h(x))=0,
$$
for any function $h$. As in the proof of Theorem \ref{theonoiseless}
we split $\Omega$ into $\Omega_+$ and $\Omega_-$ and find that
$$
\E(\|f-\t f\|^2) \leq \int_{\Omega_+} \|f-w\|^2d\rho_X^n+
8M^2n^{-r},
$$
where $w$ now stands for the solution to the least squares problem
with noisy data $(y_1,\dots,y_n)$. With the same definition of
$g=f-P_mf$, we can write
$$
f-w=g-P_m^n g- \widetilde{w},
$$
where $\widetilde{w}$ stands for the solution to the least squares problem
for the noise data $(\eta_1,\dots,\eta_n)$. Therefore
$$
\|f-P_m^nf\|^2 =\|g\|^2+\|P_m^n g+ \widetilde{w}\|^2 \leq \|g\|^2+2\|P_m^n
g\|^2+2\|\widetilde{w}\|^2 =|g\|^2+2\sum_{j=1}^m |a_j|^2 +2\sum_{j=1}^m
|d_j|^2,
$$
where $\ba=(a_j)_{j=1,\dots,m}$ is as in the proof of Theorem \ref{theonoiseless}
and $\bd=(d_j)_{j=1,\dots,m}$ is solution of the system
$$
\bG\bd=\bn, \;\;
$$
with $\bn:=(\frac 1 n\sum_{i=1}^n\eta_i L_k(x_i))_{k=1,\dots,m}=(n_k)_{k=1,\dots,m}$.
By the same arguments as in the proof of Theorem \ref{theonoiseless}, we
thus obtain
$$
\E(\|f-\t f\|^2)\leq  (1+2\e(n))e_m(f)^2+8L^2n^{-r}+ 8\sum_{k=1}^m \E(|n_k|^2).
$$
We are left to show that $\sum_{k=1}^m \E(|n_k|^2)\leq \frac {\sigma^2 m}{n}$.
For this we simply write that
$$
\E(|n_k|^2)=\frac 1 {n^2}\sum_{i=1}^n\sum_{j=1}^n \E(\eta_iL_k(x_i)\eta_j L_k(x_j)).
$$
For $i\neq j$, we have
$$
\E(\eta_iL_k(x_i)\eta_j L_k(x_j))=(\E(\eta L_k(x)))^2=0.
$$
For $i=j$, we have
\begin{align*}
\E(|\eta_iL_k(x_i)|^2)&=\E(|\eta L_k(x)|^2) \\
&=\int_X \E(|\eta L_k(x)|^2 | x) d\rho_X\\
&=\int_X\E(|\eta|^2 |x) |L_k(x)|^2d\rho_X\\
&\leq \sigma^2 \int_X |L_k(x)|^2d\rho_X=\sigma^2.
\end{align*}
It follows that $\E(|n_k|^2)\leq \frac {\sigma^2} n$, which concludes the proof. \hfill $\Box$

\section{Examples and numerical illustrations}

We now give several examples of approximation schemes for which
one can compute the quantity $K(m)$ and therefore
estimate the range of $m$ such that the condition \iref{condm} holds.
For each of these examples, we
also exhibit a deterministic sampling $(x_1,\dots,x_n)$
for which the stability property
$$
\interleave\bG-\bI\interleave \leq \frac12,
$$
or equivalently
$$
\left |\|v\|_n^2-\|v\|^2\right | \leq  \frac 1 2 \|v\|^2,\;\; v\in V_m,
$$
is ensured for the same range of $m$ (actually slightly better
by a  logarithmic factor). For the sake of simplicity, we work in the one dimensional setting,
with $X$ a bounded interval.
\nl
\nl
{\bf Piecewise constant functions.} Here $X=[a,b]$ and $V_m$ is the space
of piecewise constant functions over a partition of $X$ into
intervals $I_1,\dots,I_m$. In such a case, an orthonormal basis with respect to $L^2(X,\rho_X)$ is
given by the characteristic functions $L_k:= (\rho_X(I_k))^{-1/2}\Chi_{I_k}$,
and therefore
$$
K(m)=\max_{k=1,\ldots,m} (\rho_X(I_k))^{-1}.
$$
Given a measure $\rho_X$, the partition that minimizes $K(m)$, and
therefore allows us to fulfill \iref{condm} for the largest range of
$m$, is one that evenly distributes the measure $\rho_X$. With such
partitions, $K(m)$ reaches its minimal value
$$
K(m)=m,
$$
and \iref{condm} can be achieved with $m\sim \frac n {\log n}$.
\nl
If we now choose $n=m$ deterministic points $x_1,\ldots,x_m$ with
$x_k\in I_k$, we clearly have
$$
\|v\|_n^2=\|v\|^2,\;\; v\in V_m.
$$
Therefore the stability of the least squares problem
can be ensured with $m$ up to the value $n$ using a deterministic sample.
\nl
\nl
{\bf Trigonometric polynomials and uniform measure.} Without loss of generality, we take $X=[-\pi,\pi]$,
and we consider for odd $m=2p+1$
the space $V_{m}$ of trigonometric polynomials of degree $p$, which is spanned by the functions
$L_k(x)=e^{ikx}$ for $k=-p,\dots,p$. Assuming that $\rho_X$
is the uniform measure, this is an orthonormal basis
with respect to $L^2(X,\rho_X)$. In this example,
we again obtain the minimal value
$$
K(m)=m.
$$
Therefore \iref{condm} can be achieved with $m\sim \frac n {\log n}$.
\nl
We now consider the deterministic uniform sampling
$x_i:=-\pi+\frac {2\pi  i} n$ for $i=1,\dots,n$. With such a sampling,
one has the identity
$$
 \int_{-\pi}^\pi v(x) d\rho_X=\frac 1 {2\pi}\int_{-\pi}^\pi v(x) dx=\frac 1 n \sum_{i=1}^n v(x_i),
$$
for all trigonometric polynomials $v$ of degree $n-1$ (this is easily
seen by checking the identity on every basis element). When
$v\in V_m$ with $m=2p+1$, we know that $|v|^2$ is a trigonometric
polynomial of degree $2p$. We thus find that
$$
\|v\|_n^2=\|v\|^2,\;\; v\in V_m,
$$
provided that $2p\leq n-1$, or equivalently $m\leq n$. Therefore the
stability of the least squares problem can be ensured with $m$ up to
the value $n$ using a deterministic sample. \nl \nl {\bf Algebraic
polynomials and uniform measure.} Without loss of generality, we
take $X=[-1,1]$, and we consider $V_m=\P_{m-1}$ the space of
algebraic polynomials of degree $m-1$. When $\rho_X$ is the uniform
measure, an orthonormal basis is given by defining $L_k$ as the
Legendre polynomial of degree $k-1$ with normalization
$$
\|L_k\|_{L^\infty([-1,1])}=|L_k(1)|=\sqrt {2k-1},
$$
and thus
$$
K(m)=\sum_{k=1}^m (2k-1)=m^2.
$$
Therefore \iref{condm} can be achieved with $m\sim \sqrt {\frac n
{\log n}}$ which is a lower range compared to the previous examples.
\nl We now consider the deterministic sampling obtained by
partitioning $X$ into $n$ intervals $(I_1,\dots,I_n)$ of equal
length $\frac 2n$, and picking one point $x_i$ in each $I_i$. For
any $v\in V_m$, we may write
\begin{align*}
\left |\int_{I_i} |v(x)|^2d\rho_X - \frac 1 n |v(x_i)|^2\right | &=\left |\frac 1 2\int_{I_i} |v(x)|^2 dx - \frac 1 n |v(x_i)|^2\right | \\
& =\left |\frac 1 2\int_{I_i} (|v(x)|^2-|v(x_i)|^2) dx\right |\\
& \leq \frac 1 2\int_{I_i} \left | |v(x)|^2-|v(x_i)|^2\right | dx\\
& \leq \frac 1 2 |I_i| \int_{I_i} |(v^2)'(x)|dx \\
&= \frac 2 {n}\int_{I_i} |v'(x)v(x)| d\rho_X.
\end{align*}
Summing over $i$, it follows that
$$
\left |\|v\|_n^2-\|v\|^2\right | \leq \frac  2{n}\int_{X} |v'(x)v(x)| d\rho_X
\leq \frac 2 {n}\|v'\| \, \|v\| \leq \frac {2(m-1)^2}{n}\|v\|^2,
$$
where we have used the Cauchy-Schwarz and Markov inequalities. Therefore
the stability of the least squares problem can be ensured with $m$
up to the value $\frac{\sqrt n} 2+1$ using a deterministic sample.
\nl \nl {\bf Algebraic polynomials and Chebyshev measure.} Consider
again algebraic polynomials of degree $m-1$ on $X=[-1,1]$, now
equipped with the measure
$$
d\rho_X=\frac {dx}{\pi\sqrt{1-x^2}}.
$$
Then an orthonormal basis is given
by defining $L_k$ as the Chebyshev polynomial of degree $k-1$, with $L_1=1$ and
$$
L_k(x)=\sqrt 2\cos((k-1) \arccos x),
$$
for $k>1$, and thus
$$
K(m)=2m-1.
$$
Therefore \iref{condm} can be achieved with $m\sim \frac n {\log
n}$, which expresses the fact that least squares approximations are
stable for higher polynomial degrees when working with the
Chebyshev measure rather than with the uniform measure. \nl We now
consider the deterministic sampling obtained by partitioning $X$
into $n$ intervals $(I_1,\dots,I_n)$ of equal Chebyshev measure
$\rho_X(I_i)=\frac 1n$, and picking one point $x_i$ in each $I_i$.
For any $v\in V_m$, we may write
\begin{align*}
\left |\int_{I_i} |v(x)|^2d\rho_X - \frac 1 n |v(x_i)|^2\right | & =\left |\int_{I_i} (|v(x)|^2-|v(x_i)|^2) d\rho_X\right |\\
& \leq \int_{I_i} \left | |v(x)|^2-|v(x_i)|^2\right | d\rho_X\\
& \leq \rho_X(I_i) \int_{I_i} |(v^2)'(x)|dx \\
&= \frac 1 {n}\int_{I_i} |v'(x)v(x)| dx.\\
\end{align*}
Summing over $i$, it follows that
$$
\left |\|v\|_n^2-\|v\|^2\right | \leq \frac  1{n}\int_{X} |v'(x)v(x)| dx
\leq \frac 1 {n} \|v\| \; \(\pi \int_X |v'(x)|^2 \sqrt{1-x^2} dx\)^{1/2}.
$$
Using the change of variable $x=\cos t$, it is easily seen that the inverse estimate
$$
\int_X |v'(x)|^2 \sqrt{1-x^2} dx \leq (m-1)^2\int_X |v(x)|^2 \frac 1 { \sqrt{1-x^2}} dx,
$$
holds for any $v\in V_m$. Therefore
$$
\left |\|v\|_n^2-\|v\|^2\right | \leq \frac  1{n}\int_{X}
|v'(x)v(x)| dx \leq \frac {\pi(m-1)} {n} \|v\|^2
$$
which shows that the stability of the least squares problem can be
ensured with $m$ up to the value $\frac{ n} {2\pi}+1$ using a
deterministic sample.
\nl

Let us observe that in several practical
scenarios, the measure $\rho_X$ of the observations
may be unknown to us, therefore raising the question
of the behavior of $K(m)$ for an arbitrary measure.

It is not too difficult to check
that when the space $V_m$ is
not the trivial space of constant functions
(which is the case as soon as $m\geq 2$)
the quantity $K(m)$ may become arbitrarily large
for certain measures $\rho_X$. We leave the proof
of this general fact as an exercise for the reader,
and rather provide a simple illustration:
consider the space $V_2$ of polynomials of degree $1$
on $[-1,1]$ and the measure $\rho_X=\frac 1 {2 \e} \Chi_{[-\e,\e]}(x)dx$
where $\e>0$ is small. Then an orthonormal basis is
provided by the functions $L_0(x)=1$ and $L_1(x)=\frac {\sqrt 3} \e x$,
so that $K(m) \sim \e^{-2}$. An interesting problem is
to understand if for certain families of space $(V_m)$,
the quantity $K(m)$ can be controlled
under fairly general assumptions on the measure $\rho_X$. One typical such
assumption is the
the strong density assumption, which states that
\be
\rho_X(E) \sim |E|, \;\; E \; {\rm measurable},
\ee
where $|\cdot|$ is the Lebesgue measure.
In the case of piecewise constant functions on uniform partitions,
or for more general spline functions on uniform grids,
it is not difficult to check that this assumption implies the behavior $K(m)\sim m$.
\nl
\nl
{\bf Numerical illustration.} We
conclude with a brief numerical illustration of our theoretical
results for the setting of algebraic polynomials. Specifically, we
consider the smooth function
$f_1(x)=1/(1+25x^2)$ originally
considered by Runge to illustrate the instability of polynomial
interpolation at equispaced points, and
the non-smooth function $f_2(x) = |x|$, both
restricted to the interval $[-1,1]$.

\begin{figure}[t]
   \centering
   \begin{tabular}{cc}
   \hspace{-3mm} \includegraphics[width=.45\linewidth]{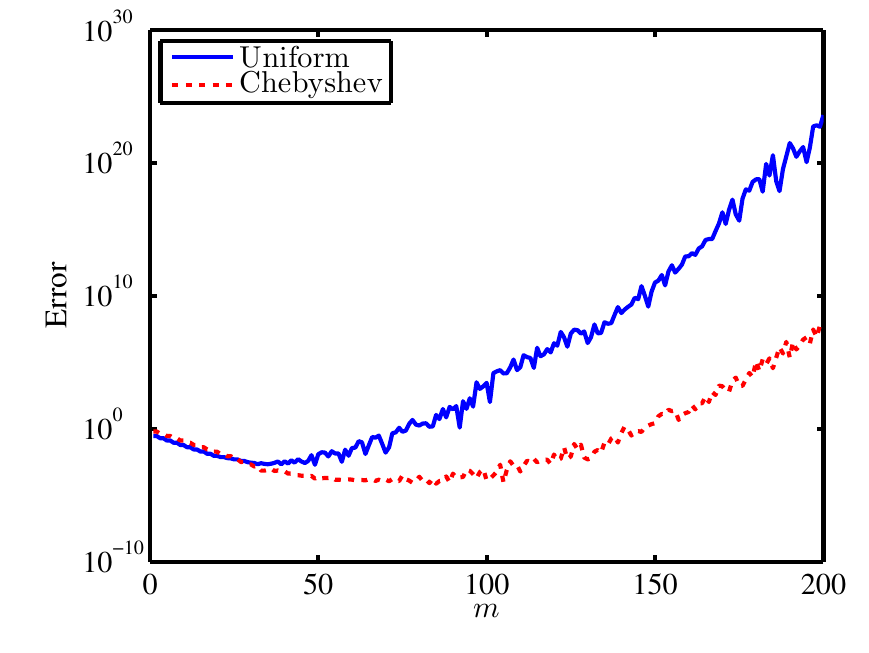} & \hspace{-5mm} \includegraphics[width=.45\linewidth]{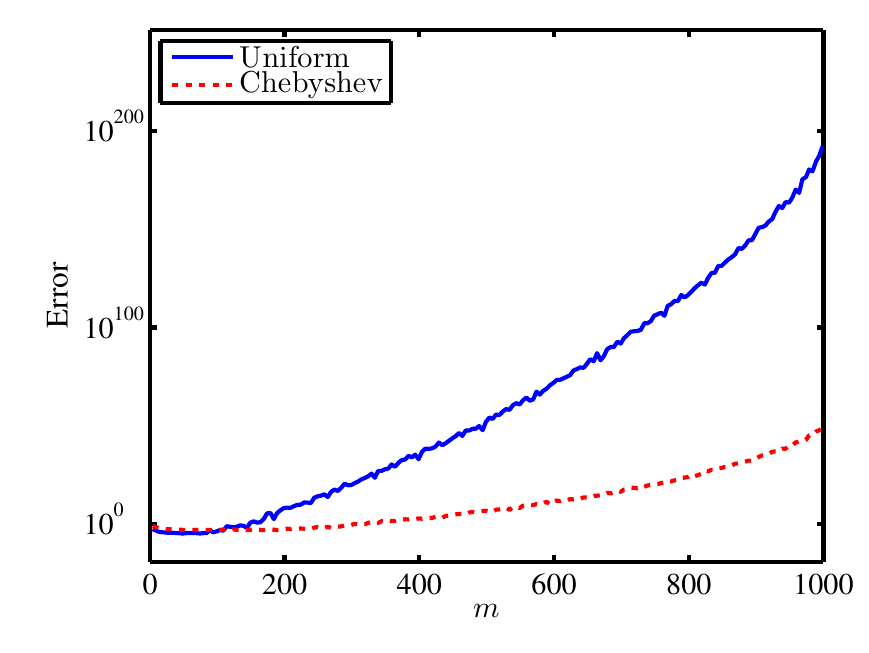} \\
   \hspace{5mm} (a) &  \hspace{2mm} (b)
   \end{tabular}
   \caption{The $L^2(X,\rho_X)$ error as $m$ varies (a) for $f_1$ and (b) for $f_2$.
   \label{fig:figure12}}
\end{figure}

For both functions, we take
$n$ i.i.d.\ samples $x_1, \ldots ,x_n$ with respect to a measure $\rho_X$ on $X=[-1,1]$ and compute the
noise-free observations $y_i=f(x_i)$. We consider either the uniform
measure $\rho_X:=\frac {dx} 2$ or the Chebyshev measure $\rho_X:=\frac {dx}{\pi\sqrt{1-x^2}}$.
In both cases, we compute the
least squares approximating polynomial of degree $m$ using these
points for a range of different values of $m \leq n$.
We then numerically compute the error in the $L^2(X,\rho_X)$ norm,
with $\rho_X$ the corresponding measure in which the sample have been drawn, using
the adaptive Simpson's quadrature rule \cite{GG} implemented in Matlab.

Figure~\ref{fig:figure12} shows the results of this simulation using $n_1 = 200$ samples for estimating $f_1$ and $n_2 = 1000$ samples for estimating $f_2$.  We observe that, in all cases, as $m$ approaches $n$ the solutions become highly inaccurate due to the inherent instability of the problem. However, we can set $m$ to be much larger
before instability starts to develop when the points are drawn with respect to the Chebyshev measure, as is expected.

\begin{figure}[t]
   \centering
   \begin{tabular}{cc}
   \hspace{-3mm} \includegraphics[width=.45\linewidth]{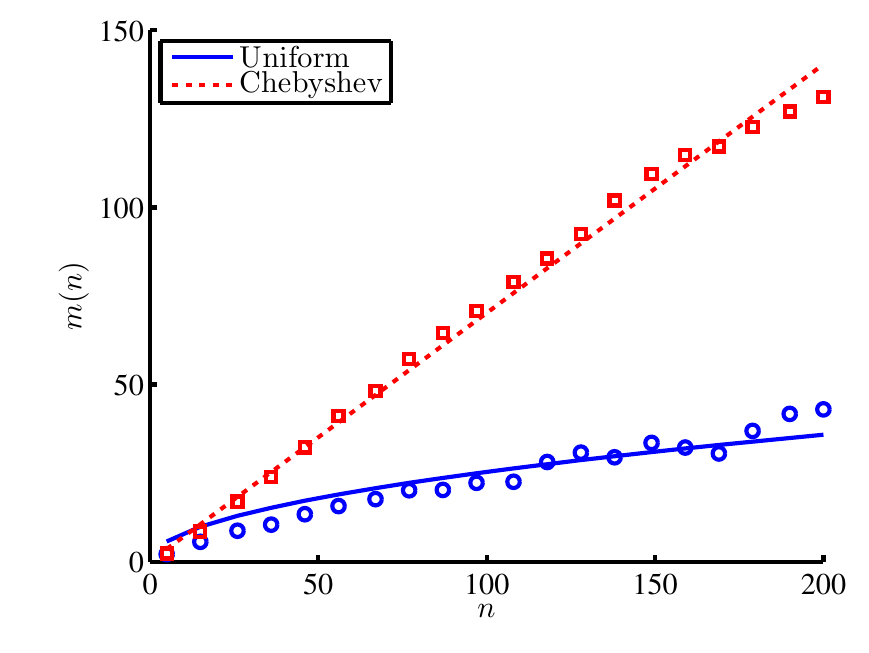} & \hspace{-5mm} \includegraphics[width=.45\linewidth]{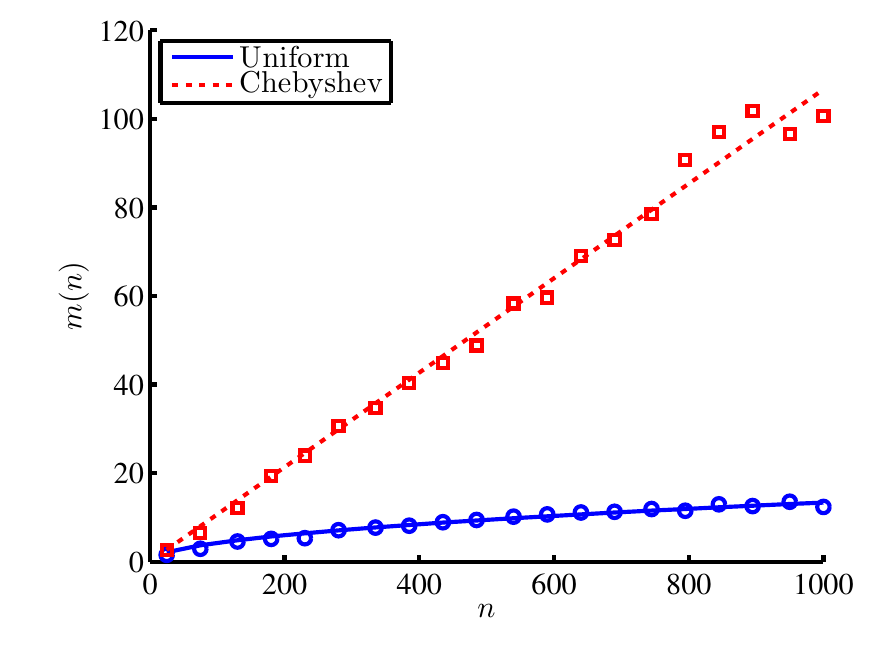} \\
   \hspace{5mm} (a) &  \hspace{2mm} (b)
   \end{tabular}
   \caption{Optimal values $m(n)$ as $n$ varies (a) for $f_1$ (comparison with
 $0.7n$ and $2.5\sqrt n$) and (b) for $f_2$
(comparison with  $0.1n$ and $0.4\sqrt n$).
   \label{fig:figure13}}
\end{figure}

Next we consider the effect of $n$ on the best choice of $m$.  Specifically, for any given sample of points we can compute the value $m(n)$ that corresponds to the polynomial degree for which
we obtain the best approximation to $f_1$ or $f_2$ and examine how this behaves as a function of $n$.
This is shown in Figure~\ref{fig:figure13}, that displays as a function of $n$ the average value of $m(n)$
over $50$ realizations of the sample, for both measures and both functions $f_1$ and $f_2$ (the
averaging has the effect of reducing oscillation in the curve $n\mapsto m(n)$ making it more readable).
We vary the sample size from
$n=1$ to $1000$ for $f_2$, but only
from $n=1$ to $200$ for the smooth function $f_1$,
since in that case the $L^2(X,\rho_X)$ error
drops below machine precision for larger values of $n$ with $m$ in the regime where the least squares problem
is stable and therefore the minimal value $m(n)$ cannot be precisely located.

We observe that, in accordance with our theoretical results, $m(n)$ behaves like $\sqrt{n}$
when the points are drawn with respect to the uniform measure,
while it behaves almost linear in $n$ when the points are drawn with respect to the Chebyshev measure.

$\;$
\nl
The authors would like to thank Lukas Meier for bringing a small error in the original proof of Theorem 1 to our attention.
\nl
\nl
Albert Cohen
\nl
Laboratoire Jacques-Louis Lions
\nl
Universit\'e Pierre et Marie Curie
\nl
4, Place Jussieu, 75005 Paris, France
\nl
cohen@ann.jussieu.fr
\nl
\nl
Mark Davenport
\nl
School of Electrical and Computer Engineering
\nl
Georgia Institute of Technology
\nl
777 Atlantic Drive NW
\nl
Atlanta, GA 30332, USA
\nl
mdav@gatech.edu
\nl
\nl
Dany Leviatan
\nl
Raymond and Beverly Sackler School of Mathematics
\nl
Tel Aviv University
\nl
69978, Tel Aviv, Israel
\nl
leviatan@post.tau.ac.il
\nl

\end{document}